\documentclass{amsart}
\usepackage[arrow,matrix]{xy}

 \newcommand{\R}{\mathbb R}
 
 \DeclareMathOperator{\Int}{Int}
\DeclareMathOperator{\cl}{cl} \DeclareMathOperator{\tr}{tr}

\newcommand{\nphi}{\tilde\varphi}

\theoremstyle{plain} \newtheorem{thm}{Theorem}
\newtheorem{cor}[thm]{Corollary} 
\newtheorem{lemma}[thm]{Lemma}

\theoremstyle{definition}

\theoremstyle{remark} \newtheorem{remark}[thm]{Remark}

\begin{document}
    
\title[A fixed point theorem for bounded dynamical systems]{A fixed
point theorem for bounded dynamical systems}

\author{David Richeson} \author{Jim Wiseman} \address{Dickinson
College\\Carlisle, PA 17013} \email{richesod@dickinson.edu}
\address{Northwestern University\\Evanston, IL 60208}
\email{jimw@math.northwestern.edu}

\begin{abstract}
    We show that a continuous map or a continuous flow on $\R^{n}$
    with a certain recurrence relation must have a fixed point. 
    Specifically, if there is a compact set $W$ with the property that
    the forward orbit of every point in $\R^{n}$ intersects $W$ then
    there is a fixed point in $W$.  Consequently, if the omega limit
    set of every point is nonempty and uniformly bounded then there is
    a fixed point.
\end{abstract}

\maketitle

In this note we will prove a fixed point theorem that holds for both
discrete dynamical systems ($f$ is a continuous map) and continuous
dynamical systems ($\varphi^{t}$ is a continuous flow).  We will show
that if every point in $\R^{n}$ returns to a compact set, then there
must be a fixed point.  This investigation began in an attempt to
answer a related question about smooth flows posed by Richard
Schwartz.

We will prove the theorem for maps, and derive the theorem for flows
as a consequence.  In fact, many of the definitions and proofs follow
analogously for both cases.  When that is the case we will just refer
to the ``dynamical system,'' with the recognition that the statement
applies to both flows and maps.  Where necessary, separate definitions
and proofs will be included.

We now introduce the notion of a window; a compact set $W$ is a {\em
window} for a dynamical system on $X$ if the forward orbit of every
point $x\in X$ intersects $W$.  If a dynamical system has a window
then we will say that it is {\em bounded}.

We will prove the following fixed point theorem.

\begin{thm} \label{thm:fpt}
    Every bounded dynamical system on $\R^n$ has a fixed point.
\end{thm}

The following corollaries are elementary applications of Theorem
\ref{thm:fpt}.

\begin{cor} \label{cor:inW}
    If $W$ is a window for a dynamical system on $\R^{n}$, then there
    is a fixed point in $W$.
\end{cor}

\begin{cor}
    If there is a compact set $K$ such that
    $\emptyset\ne\omega(x)\subset K$ for all $x\in\R^{n}$, then the
    dynamical system has a fixed point.
\end{cor}

We will need the following theorem.

\begin{thm}[Lefschetz Fixed Point Theorem] \label{lef-thm}
    Let $f:M\to M$ be a continuous map of an $n$-dimensional manifold
    (with or without boundary), and let $f_{k}:H_{k}(M;\R)\to
    H_{k}(M;\R)$ be the induced map on homology.  If
    $\sum_{k=0}^{n}(-1)^{k}\tr(f_{k}) \neq 0$, then $f$ has a fixed
    point.
\end{thm}

The main result in this paper concerns dynamical systems on $\R^{n}$,
but some of the results hold more generally.  Unless otherwise stated,
our dynamical system is defined on a locally compact topological
space.

We have the following result.  The proof is similar to one by Conley
\cite[\S~II.5]{C}.

\begin{lemma} \label{lem:attr}
    Every bounded dynamical system has a forward invariant window.
\end{lemma}

\begin{proof}
    Suppose our dynamical system is a continuous map $f$.  Because $f$
    is bounded, there is a window $W$.  Without loss of generality, we
    may assume that the forward orbit of each point visits the
    interior of $W$.  If not, then we may replace $W$ by any compact
    neighborhood of $W$.  For each $x\in W$ there exists an $n_x> 0$
    for which $f^{n_x}(x)\in\Int W$.  Clearly there is an open
    neighborhood $U_x$ of $x$ such that $f^{n_x}(y)\in\Int W$ for all
    $y\in U_x$.  The sets $\{U_x:x\in W\}$ form an open cover of $W$. 
    Since $W$ is compact there is a finite subcover,
    $\{U_{x_1},\ldots,U_{x_m}\}$.  Let
    $n=\max\{n_{x_k}:k=1,\ldots,m\}$.  It follows that
    $W_0=\bigcup_{k=0}^{n}f^k(W)$ is a forward invariant window.
    
    The proof for a flow $\varphi^{t}$ is similar.
\end{proof}

\begin{proof}[Proof of Theorem \ref{thm:fpt} for maps]
    Let $f:\R^n\to\R^n$ be a bounded, continuous map with a forward
    invariant window $W_{0}$.  Let $B$ be a closed ball containing
    $W_0$ in its interior.  We begin by constructing a manifold with
    boundary $N$ containing $B$ such that $f(N)\subset\Int N$.
    
    Arguing as in the proof of Lemma~\ref{lem:attr}, there is a
    positive integer $n$ such that for each $x$ in $B$, the set $x\cup
    f(x) \cup \dots \cup f^n(x)$ intersects $W_0$.  Since $W_0$ is
    forward invariant, it follows that $f^n(B) \subset W_0$.  The set
    $\bigcup_{k=0}^{n}f^{k}(B)$ is forward invariant; following
    \cite[Thm.  3.3(a)]{A} we will enlarge this set slightly so that
    it maps into its interior.  Define the set-valued map $V_r:\R^n
    \to \R^n$ sending a point $x$ to the closed ball of radius $r$
    centered at $x$.  Since $f^n(B) \subset \Int B$, Lemma 3.2 of
    \cite{A} tells us that for $\delta > 0$ sufficiently small
    $(V_{\delta}\circ f)^n(B) \subset \Int B$ (the lemma is actually
    stated for compact spaces, but the same proof remains valid in
    this case).  Thus, the set $B_0 = B\cup (V_{\delta}\circ f)(B)
    \cup \dots \cup (V_{\delta}\circ f)^{n-1}(B)$ has the property
    $f(B_{0})\subset\Int B_{0}$.  Finally, there exists a submanifold
    with boundary $N$ sufficiently close (in the Hausdorff topology)
    to $B_0$ containing $B$ such that $f(N) \subset \Int N$.
    
    We claim that $f_{k}:H_{k}(N;\R)\to H_{k}(N;\R)$ is nilpotent for
    $k\ne 0$.  Again, there is a positive integer $m$ such that
    $f^m(N) \subset W_0 \subset B$.  Thus the map $f^m: N \to N$
    factors through $B$, i.e., the following diagram commutes: \\

    
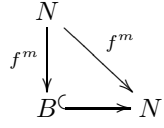
\hspace{1.5in} \xymatrix{ N \ar[dr]^{f^m} \ar[d]_{f^m} \\
    B \ar@{^{(}->}[r] & N}\\
    Therefore the map $(f^m)_{\ast} : H_{\ast}(N;\R) \to
    H_{\ast}(N;\R)$ factors through $H_{\ast}(B;\R)$.  Since $B$ is
    contractible, its homology consists of an $\R$ in dimension zero
    and zeroes elsewhere, so $(f^m)_{\ast}$ is the identity in
    dimension zero and the zero map elsewhere.  That is, for $k\neq
    0$, $0 = (f^m)_k = (f_k)^m$, i.e., $f_k$ is nilpotent.  Since the
    trace of a nilpotent matrix is zero, the alternating sum of the
    traces of $f_\ast$ is 1.  Thus, $N$ must contain a fixed point.
\end{proof}

\begin{proof}[Proof of Theorem \ref{thm:fpt} for flows]
    For $s\in [0,1]$, define $f_s: \R^n\to \R^n$ to be the time-$s$
    map of $\varphi^t$ (i.e., $f_s(x) = \varphi^s(x)$).  If $W$ is a
    window for the flow, then the set $W' = \bigcup_{t=0}^1
    \varphi^t(W)$ is a window for each $f_s$.  Thus
    Corollary~\ref{cor:inW} for maps tells us that every $f_s$ has a
    fixed point in $W'$.

    Take a sequence of $s$'s tending to zero and let $x$ be a limit
    point of the corresponding sequence of fixed points of the
    $f_s$'s; we claim that $x$ is a fixed point of $\varphi^t$.

    If not, then there exist a neighborhood $U$ of $x$ and a time
    $t_0$ such that $\varphi^{t_0}(U) \cap U = \emptyset$.  But every
    neighborhood of $x$ must contain periodic orbits, which is a
    contradiction.
\end{proof}

\begin{remark}
    The same proof can apply in other spaces.  What we need is that
    the forward invariant window $W_0$ have a compact acyclic
    neighborhood $B$ (i.e., $\tilde H_*(B;\R) = 0$).  In the case of a
    nonsingular flow on $S^1$, for example, $W_0 = S^1$ and the proof
    fails.
\end{remark}

\begin{remark}
    There is another proof of Theorem~\ref{thm:fpt} in the case of a
    smooth flow $\varphi^t$, the outline of which we will present now. 
    We may add a point $p$ at infinity to $\R^n$ to obtain $S^n$;
    $\varphi^t$ induces a smooth flow $\nphi^t$ on $S^n$ with a fixed
    point at $p$.  Because there exists a smooth Lyapunov function for
    $\nphi^t$, we can find a submanifold with boundary $L$ such that
    $\nphi^t(L) \subset \Int L$ for all $t\geq 0$ and $\bigcap_{t\geq
    0}\nphi^{-t}(\cl(S^n\backslash L)) = \{p\}$.  The set
    $\cl(S^n\backslash L)$ is contractible (the homotopy is given by
    running the flow backwards, so that every point goes to $p$). 
    Thus $\cl(S^n\backslash L)$ is diffeomorphic to a closed $n$-disk;
    by the Schoenflies theorem, the same is true of $L$.  Therefore
    $L$ is acyclic, and we can proceed as above.
\end{remark}

A closed subset $C\subset X$ is an {\em attracting neighborhood} for
$f$ if $f(C)\subset \Int C$.  A set $A$ is an {\em attractor} provided
there is an attracting neighborhood $C$ such that $A =
\bigcap_{k\geq0}f^k(C)$.  (See \cite{C}, \cite[Ch.  3]{A} for details. 
Attractors for flows are defined similarly.)  An attractor $A$ is a
{\em global attractor} if $\emptyset\ne\omega(x)\subset A$ for all
$x\in X$.

\begin{cor}
    Every bounded dynamical system has a unique global attractor, $A$. 
    If the dynamical system is a homeomorphism or a flow on $\R^n$
    then the reduced \v{C}ech cohomology of $A$,
    $\Tilde{\Check{H}}^{\ast}(A)$, is zero.
\end{cor}

\begin{proof}
    For a map $f$, the set $A=\bigcap_{n\ge 0}f^{n}(W_{0})$ is clearly
    a global attractor.  Suppose $A^{\prime}$ is another global
    attractor and $N$ and $M$ are attracting neighborhoods for $A$ and
    $A^{\prime}$ respectively.  Since $N$ and $M$ are both windows,
    there exists $n\ge 0$ such that $f^n(N) \subset M$ and $f^n(M)
    \subset N$.  Thus, $A=A^{\prime}$.  We argue analogously for a
    flow $\varphi^t$.
    
    If the dynamical system is invertible, then $A$ is equal to the
    intersection of a nested sequence of open, contractible sets.  For
    a map $f$, $A = \bigcap_{k\geq 0}f^{km}(\Int B)$ with $m,B$ as
    above.  In the flow case, for any ball $B^{\prime}$ containing
    $W_{0}$ in its interior, there is a $t_{0}>0$ such that
    $\varphi^{t_{0}}(B^{\prime})\subset W_{0}\subset\Int B^{\prime}$. 
    So, $A = \bigcap_{k\ge 0}\varphi^{kt_{0}}(\Int B^{\prime})$.  Thus
    $\Tilde{\Check{H}}^{\ast}(A) = 0$.
\end{proof}

Of course, if the dynamical system is a homeomorphism or a flow then
one can give the analogous results for $f^{-1}$ or $\varphi^{-t}$.

\begin{cor}
    Suppose $f:\R^{n}\to\R^{n}$ is a homeomorphism (respectively
    $\varphi^{t}$ is a flow).  If $f^{-1}$ is bounded (resp. 
    $\varphi^{-t}$) then there is a fixed point.  In particular, if
    there is a compact set $K$ such that
    $\emptyset\ne\alpha(x,f)\subset K$ (resp. 
    $\emptyset\ne\alpha(x,\varphi^{t})\subset K$) for all $x\in\R^{n}$
    then there is a fixed point.
\end{cor}

Finally, we have the following interesting corollary.

\begin{cor}
    Every dynamical system on a non-compact, locally compact space has
    a point whose forward orbit is not dense.
\end{cor}

\begin{proof}
    If the forward orbit of every point is dense then any compact set
    $W$ with nonempty interior is a window.  But Lemma~\ref{lem:attr}
    shows that the existence of a window implies the existence of a
    compact, forward invariant set, which cannot be.
\end{proof}

\bibliographystyle{cite} 
\bibliography{fptbib}
\end{document}